\documentclass[secthm]{paper}
\usepackage{amssymb,amsmath}
\usepackage{amsfonts}
\usepackage{latexsym,diagrams}
\input{defs_allg.tex.input}

\begin{document}
\hyphenation{co-qua-si-tri-an-gu-lar}
\begin{frontmatter}
\title{Triangular braidings and pointed Hopf algebras\thanksref{PhD}}
\thanks[PhD]{This work is part of the authors PhD thesis written under the supervision of Professor H.-J. Schneider.}
\author{Stefan Ufer\thanksref{BayStaat}}
\address{Mathematisches Institut der Universit\"at M\"unchen, Theresienstr. 39,\\ 80333 M\"unchen, Germany}
\thanks[BayStaat]{Partially supported by Graduiertenf\"orderung des bayerischen Staates and by his parents}
\begin{abstract}
We consider an interesting class of braidings defined in \cite{Ich_PBW} by a combinatorial property. We show that it consists exactly of those braidings that come from certain Yetter-Drinfeld module structures over pointed Hopf algebras with abelian coradical.\\
As a tool we define a reduced version of the FRT construction. For braidings induced by $U_q(\mathfrak{g})$-modules the reduced FRT construction is calculated explicitly.
\end{abstract}
\end{frontmatter}

\section{Introduction}

In \cite{Ich_PBW} the author generalized a PBW theorem of Kharchenko \cite{kh.b} for certain pointed Hopf algebras to a class of braided Hopf algebras. The central feature of these braided Hopf algebras is that they are generated by a subspace of primitive elements that has a so-called triangular braiding. This notion of triangularity was defined in \cite{Ich_PBW} by a combinatorial property. It seems to be the natural context for the proof of the PBW theorem given there.\par
The combinatorial description is very helpful for the proof of the PBW theorem, but does not give a conceptual understanding of these braidings. The main examples of triangular braidings, namely Yetter-Drinfeld modules over abelian groups and integrable modules over quantum enveloping algebras, come in some sense from pointed Hopf algebras with abelian coradical. As also the original PBW theorem of Kharchenko is proved for such Hopf algebras, it is an interesting question if there is a deeper connection between these Hopf algebras (including quantum groups) and triangular braidings.\par
This work shows that there is indeed a close connection. We show (Theorem \ref{thm_hopfalg}) that right triangular braidings are exactly those braidings that are induced by Yetter-Drinfeld modules over pointed Hopf algebras with abelian coradical which are completely reducible as $kG(H)$-modules.\par
It is a well-known fact \cite{FRT} that any braiding $c$ on a finite-dimensional vector space $M$ can be realized as a Yetter-Drinfeld braiding over some bialgebra $A(c)$, the FRT bialgebra. If the braiding is rigid (and we prove that triangular braidings are rigid) it can even be realized as a Yetter-Drinfeld braiding over a Hopf algebra $H(c)$. This was first proved by Lyubashenko for symmetries \cite{Lyubashenko}; for the general case see \cite{Sbg_cqt}.\par
Nevertheless for our purpose the constructions from \cite{FRT,Lyubashenko,Sbg_cqt} do not seem to be the right tool to study properties of the braiding. In \cite{Radford_pointed} Radford defines a reduced version $A^{red}$ of the FRT bialgebra and uses it as a tool to obtain information on braidings induced by Yetter-Drinfeld modules over pointed bialgebras.\par
In this paper we give a different construction for Radford's reduced FRT bialgebra. This construction leads in a natural way to a reduced FRT Hopf algebra if the braiding is rigid. We prove an analogue of Radford's central result in the case of the reduced FRT Hopf algebra. This is an important tool for the proof of Theorem \ref{thm_hopfalg}. We apply Theorem \ref{thm_hopfalg} to show that there are braidings which are left, but not right triangular.\par
In the case of braidings constructed from the quasi-\Rmatrix{} of $U_q(\mathfrak{g})$ the usual FRT constructions lead to interesting examples of Hopf algebras, namely the quantized function algebras. In the last section we determine the reduced FRT Hopf algebra for these braidings and we obtain pointed Hopf algebras which are closely related to the non-positive part $U_q^{\leq 0}(\mathfrak{g})$ of $U_q(\mathfrak{g})$.\\[0,5cm]
The paper is organized as follows: Section \ref{sect_cqtbialg} contains the definition of coquasitriangular bialgebras, their right radical and their reduced versions. Furthermore we recall some facts on Yetter-Drinfeld modules. Section \ref{sect_FRT} deals with the FRT constructions and their reduced versions. In Section \ref{sect_redpointed} results of Radford are generalized to the reduced FRT Hopf algebra. Section \ref{sect_triang} contains the central theorem about triangular braidings and finally in Section \ref{sect_uqg} the reduced FRT Hopf algebra for braidings of integrable $U_q(\mathfrak{g})$-modules is determined.\par
Throughout the paper $k$ is a field, tensor products are always taken over $k$. We use Sweedler's notation for the comultiplication $\Delta(x) = x\s1\otimes x\s2$.

\section{Coquasitriangular bialgebras and their right radical}
\label{sect_cqtbialg}
In this section we will recall some well known facts about coquasitriangular bialgebras and define the right radical of those bialgebras. If $A$ is an algebra and $C$ is a coalgebra denote the convolution product on the algebra $\Hom_k(C,A)$ by $\star$.
% **** CQT Bialgebras ****
\begin{defn}
\label{def_cqt}
A coquasitriangular bialgebra $(H,\nabla,\eta,\Delta,\eps,r)$ is a bialgebra together with a convolution invertible linear form $r\in (H\otimes H)^*$ satisfying
\[ \nabla^{op} = r\star \nabla \star r^{-1} \]
and
\[ r\circ(\nabla\otimes \id_H) = r_{13}\star r_{23}\:\:\mbox{and}\:\: r\circ(\id_H\otimes\nabla) = r_{13}\star r_{12}\]
where we define $r_{12},r_{23},r_{13}\in (H\otimes H\otimes H)^*$ by
\[r_{12} := r\otimes\eps,\:\:r_{23}:=\eps\otimes r\:\:\mbox{and}\:\: r_{13}(g\otimes h\otimes l) := \eps(h) r(g\otimes l)\]
for all $g,h,l\in H$.
\end{defn}
\begin{rem}
Let $(H,\nabla,\eta,\Delta,\eps,S,r)$ be a coquasitriangular Hopf algebra.
\begin{enumerate}
\item $S^2$ is a coinner automorphism of $H$. In particular $S$ is invertible.
\item \[r\circ(\eta\otimes \id_H) = \eps\:\:\mbox{and}\:\: r\circ(\id_H\otimes\eta)=\eps\]
\item \[r\circ(S\otimes \id_H) = r^{-1},\:\: r^{-1}\circ(\id_H\otimes S) = r,\:\: r\circ(S\otimes S) = r\]
\end{enumerate}

If $H$ is a coquasitriangular bialgebra, the second and third axiom from Definition \ref{def_cqt} read for $a,b,c\in H$:
\[ r(ab\otimes c) = r(a\otimes c\s1)r(b\otimes c\s2),\:\: r(a\otimes bc) = r(a\s2\otimes b)r(a\s1\otimes c)\]
\end{rem}

% **** The right radical ****
\begin{defn}
Let $H$ be a coquasitriangular bialgebra. Define the right radical of $H$ as
\[ J_H := \{h\in H|r(H\otimes h)=0\}.\]
\end{defn}
\begin{lem}
Let $H$ be a coquasitriangular bialgebra.
\begin{enumerate}
\item The right radical is a biideal in $H$.
\item If $H$ is a Hopf algebra, then the right radical is stable under $S$ and $S^{-1}$.
\end{enumerate}
\end{lem}
\begin{pf}
The properties of $r$ imply that the map
\[ H\rightarrow \left(H^\circ\right)^{cop}, h\mapsto r(-\otimes h)\]
is a (well-defined) morphism of bialgebras (resp. Hopf algebras). $J_H$ is the kernel of this map.
\end{pf}
We will find the following lemma usefull:
\begin{lem}
\label{lem_radical_gen}
Let $H$ be a coquasitriangular bialgebra generated as an algebra by a subset $X\subset H$. Then for every coideal $J\subset H$ we have $r(H\otimes J)=0$ if and only if $r(X\otimes J)=0$ and thus
\[ J_H = \sum\{J|J\subset H\:\mbox{is a coideal with}\:r(X\otimes J)=0\}\]
\end{lem}
\begin{pf}
The first part follows from the definition of the $r$-form. The second part is trivial.
\end{pf}

% **** Reducing cqt bialgebras ****
Now we define the reduced version of a coquasitriangular bialgebra (Hopf algebra).
\begin{defn}
Let $H$ be a coquasitriangular bialgebra (Hopf algebra) and $J_H$ its right radical. Define
\[H^{red} := H/J_H\]
the factor bialgebra (Hopf algebra).
\end{defn}
Note that if $H$ is a coquasitriangular Hopf algebra, then $H$ and $H^{red}$ have bijective antipodes.\par

% **** Yetter-Drinfeld modules over bialgebras ****
Instead of comodules over the FRT bialgebra we will consider Yetter-Drinfeld modules over the reduced FRT bialgebra. We need a more general notion of Yetter-Drinfeld modules than usual.
\begin{defn}
Let $H$ be a bialgebra. An left $H$-module and left $H$-comodule $M$ is called a (left-left) $H$ Yetter-Drinfeld module if for all $h\in H,m\in M$ we have
\[ \left(h\s1 m\right)_\sm1 h\s2 \otimes \left(h\s1 m\right)_\s0 = h\s1 m\sm1 \otimes h\s2 m\s0 \]
\end{defn}
If $H$ is a Hopf algebra this condition is equivalent to the usual Yetter-Drinfeld condition
\[ \delta(hm) = h\s1 m\sm1 S\left(h\s3\right) \otimes h\s2m\s0. \]
In any case we have a natural transformation
\[ c_{M,N}:M\otimes N\rightarrow N\otimes M,\:\: c(m\otimes n) = m\sm1 n\otimes m\s0\:\:\mbox{for}\:\:M,N\in{}^H_H\mYD\]
that satisfies the hexagon equations and the braid equation. If $H$ has a skew antipode $\bar{S}$ then ${}^H_H\mYD$ is a braided monoidal category with inverse of the braiding given by
\[ c^{-1}_{N,M}:N\otimes M\rightarrow M\otimes N, \:\: c^{-1}(n\otimes m) =  m\s0\otimes \bar{S}\left(m\sm1\right) n.\]

% **** Yetter-Drinfeld modules and comodules over cqt bialgs ****
\begin{rem}
\label{rem_categories}
Let $H$ be a coquasitriangular bialgebra. The category ${}^H\mM$ is a braided monoidal category with braiding
\[ c_{M,N}:M\otimes N\rightarrow N\otimes M,\:\: c(m\otimes n)= r\left(n\sm1\otimes m\sm1\right)n\s0\otimes m\s0.\]
Every $H$-comodule $M$ becomes a Yetter-Drinfeld module over $H$ with the action
\[\forall h\in H,m\in M:\:h\cdot m := r(m\sm1\otimes h)m\s0.\]
This defines a functor ${}^H\mM\rightarrow {}^H_H\mYD$ that is compatible with the monoidal structure and the braiding. Similarily we can define a functor ${}^H\mM\rightarrow {}^{H^{red}}_{H^{red}}\mYD$ that has the same properties.
\end{rem}
% **** FRT-construction ****
\section{FRT construction}
\label{sect_FRT}
In this section we will introduce the well known FRT constructions and reduced versions of these. For proofs the reader is referred to \cite{Sbg_cqt,Kassel}.
\begin{defn}
A \emph{braided vector space} $(M,c)$ is a vector space $M$ together with an automorphism $c$ of $M\otimes M$ that satisfies the braid equation
\[ (\id_M\otimes c)(c\otimes\id_M)(\id_M \otimes c) =
(c\otimes\id_M)(\id_M\otimes c)(c\otimes\id_M).\]
\end{defn}
\begin{thm}
\label{thm_FRTbialg}
Let $(M,c)$ be a finite dimensional braided vectorspace. 
\begin{itemize}
\item There is a coquasitriangular bialgebra $A(c)$ such that $M$ is a left $A(c)$-comodule and the braiding $c$ equals the braiding on $M$ induced by the coquasitriangular structure of $A(c)$.
%\item For all coquasitriangular bialgebras $B$ having $M$ as a left comodule such that the braiding induced by the coquasitriangular structure equals $c$ there is a unique morphism of coquasitriangular bialgebras $\phi:A(c)\rightarrow B$ such that
%\[ \delta_B = (\phi \otimes \id_M)\delta_{A(c)}.\]
\item For all bialgebras $B$ having $M$ as a Yetter-Drinfeld module such that the induced braiding equals $c$ there is a unique morphism of bialgebras $\phi:A(c)\rightarrow B$ such that
\[ \delta_B = (\phi \otimes \id_M)\delta_{A(c)}\:\:\mbox{and}\:\:\forall u\in A(c),m\in M: u\cdot m= \phi(u)\cdot m.\]
The $A(c)$-action on $M$ is defined as in Remark \ref{rem_categories}.
\end{itemize}

The algebra $A(c)$ is generated by the smallest subcoalgebra $C\subset A(c)$ satisfying $\delta_{A(c)}(M)\subset C\otimes M$.
\end{thm}
\begin{pf}
The proof is similar to the one for Theorem \ref{thm_FRTHA}.
\end{pf}

\begin{defn}
\label{defn_evdb}
For any finite dimensional vectorspace $M$ define $k$-linear maps
\[ \ev: M^*\otimes M \rightarrow k,\: \ev(\phi\otimes m) := \phi(m)\]
\[ \db: k\rightarrow M\otimes M^*,\: \db(1) := \sum_{i=1}^n m_i\otimes m^i\]
where $m_1,\ldots,m_n$ form a basis of $M$ and $m^1,\ldots,m^n$ is the dual basis of $M^*$.
\end{defn}

\begin{defn}
\label{defn_rigid}
A braided vectorspace $(M,c)$ will be called \emph{rigid} if it is finite dimensional and the map $c^\flat:M^*\otimes M\rightarrow M\otimes M^*$ defined by
\[ c^\flat := (\ev\otimes\id_{M\otimes M^*})(\id_{M^*}\otimes c\otimes \id_{M^*})(\id_{M^*\otimes M}\otimes \db)\]
is an isomorphism.
\end{defn}

Now we can formulate the Hopf algebra version of Theorem \ref{thm_FRTbialg}
\begin{thm}
\label{thm_FRTHA}
Let $(M,c)$ be a rigid braided vectorspace. 
\begin{itemize}
\item There is a coquasitriangular Hopf algebra $H(c)$ such that $M$ is a left $H(c)$-comodule and the braiding $c$ equals the braiding on $M$ induced by the coquasitriangular structure of $H(c)$.
%\item For all coquasitriangular Hopf algebras $H$ having $M$ as a left comodule such that the braiding induced by the coquasitriangular structure equals $c$ there is a unique morphism of coquasitriangular Hopf algebras $\phi:H(c)\rightarrow H$ such that
%\[ \delta_H = (\phi \otimes \id_M)\delta_{H(c)}.\]
\item For all Hopf algebras $H$ having $M$ as a Yetter-Drinfeld module such that the induced braiding equals $c$ there is a unique morphism of Hopf algebras $\psi:H(c)\rightarrow H$ such that
\[ \delta_H = (\psi \otimes \id_M)\delta_{H(c)}\:\:\mbox{and}\:\:\forall u\in H(c), m\in M: u\cdot m = \psi(u)\cdot m\]
\end{itemize}
Let $C\subset H(c)$ be the smallest subcoalgebra satisfying $\delta_{H(c)}(M)\subset C\otimes M$. Then the algebra $H(c)$ is generated by $C+S(C)$.
\end{thm}

\begin{pf}
The existence of the coquasitriangular Hopf algebra $H(c)$ is proved in \cite[Theorem  3.2.9]{Sbg_cqt}. The universal property we give is bit stronger than the one given there, but using the characterization of $H(c)$ by generators and relations in \cite[Lemma 3.2.11]{Sbg_cqt} it is easy to check our stronger version.\qed
\end{pf}

% **** Reduced FRT Constructions ****
For a braided vectorspace $(M,c)$ we define the reduced FRT bialgebra by
\[ A^{red}(c) := (A(c))^{red}\]
and if $(M,c)$ is rigid define the reduced FRT Hopf algebra by
\[ H^{red}(c) := (H(c))^{red}.\]

\begin{defn}
Let $H$ be a bialgebra and $M_1,\ldots,M_s$ $H$-modules. We will call $H$ $M_1,\ldots,M_s$\emph{-reduced} if $(0)$ is the only coideal of $H$ annihilating all the $M_i$.
\end{defn}

The reduced FRT constructions are characterized by universal properties:

\begin{thm}
\label{thm_redFRTuniv}
Let $(M,c)$ be a finite dimensional braided vectorspace.
\begin{enumerate}
\item $M$ is a Yetter-Drinfeld module over $A^{red}(c)$ such that the induced braiding is $c$. $A^{red}(c)$ is $M$-reduced.
\item For every bialgebra $A$ having $M$ as a Yetter-Drinfeld module such that the induced braiding is $c$ and such that $A$ is $M$-reduced there is a unique monomorphism of bialgebras $\phi:A^{red}(c)\rightarrow A$ such that
\[ \delta_A = (\phi\otimes M)\delta_{A^{red}(c)} \:\:\mbox{and}\:\:\forall u\in A^{red}(c),m\in M: u\cdot m = \phi(u)\cdot m.\]
\item Assume $(M,c)$ is rigid. $M$ is a Yetter-Drinfeld module over $H^{red}(c)$ such that the induced braiding is $c$. $H^{red}(c)$ is $M,M^*$-reduced.
\item Assume $(M,c)$ is rigid. For every Hopf algebra $H$ having $M$ as a Yetter-Drinfeld module such that the induced braiding is $c$ and such that $H$ is $M,M^*$-reduced there is a unique monomorphism of Hopf algebras $\psi:H^{red}(c)\rightarrow H$ such that
\[ \delta_H = (\psi\otimes M)\delta_{H^{red}(c)} \:\:\mbox{and}\:\:\forall u\in H^{red}(c),m\in M: u\cdot m = \psi(u)\cdot m\]
\end{enumerate}
\end{thm}
\begin{pf}
The right radical $J_{A(c)}$ (resp. $J_{H(c)}$) is the maximal coideal annihilating $M$ (resp. $M$ and $M^*$): If $J$ is another coideal annihilating $M$ (resp. $M$ and $M^*$) then Lemma \ref{lem_radical_gen} and Theorems \ref{thm_FRTbialg}, \ref{thm_FRTHA} imply that $r(A(c)\otimes J)=0$ (resp. $r(H(c)\otimes J)=0$) and thus $J$ is contained in the right radical. Parts one and three follow.
Now we deal with parts two and four simultaneously:
Using the universal property of the FRT constructions we find morphisms of bialgebras
\[ \hat{\phi}:A(c)\rightarrow A\:\:\mbox{and}\:\:\hat{\psi}:H(c)\rightarrow H.\]
which are compatible with the action and the coaction. The right radical of $A(c)$ (resp. $H(c)$) is the maximal coideal annihilating $M$ (resp. $M$ and $M^*$) by Lemma \ref{lem_radical_gen}. Thus the image of the right radical under $\hat{\phi}$ (resp. $\hat{\psi}$) is again a coideal annihilating $M$ (resp. $M,M^*$). As $A$ (resp. $H$) is $M$ (resp. $M,M^*$)-reduced we see that the right radical is mapped to $(0)$. This means that $\hat{\phi}$ (resp. $\hat{\psi}$) factorize over the reduced FRT-constructions. These induced maps are compatible with action and coaction. Injectivity of the induced maps follows because of the maximality of the right radical mentioned above.
\qed\end{pf}

\begin{rem}
In \cite{Radford_pointed} Radford defines a reduced FRT-bialgebra $A^{red}(R)$ for Yang-Baxter operators $R$ on finite dimensional vectorspaces $M$, that is automorphisms $R$ of $M\otimes M$ satisfying the quantum Yang-Baxter equation
\[ R_{12}R_{13}R_{23} = R_{23}R_{13}R_{12}.\]
It is well known that $R$ satisfies the quantum Yang-Baxter equation if and only if $c:= R\tau$ is a braiding (satisfies the braid equation). It is easy to see from the universal properties of $A^{red}(c)$ and $A^{red}(R)$ that if $c = R\tau$ we have
\[ A^{red}(c) \isom A^{red}(R)^{cop}. \]
\end{rem}

\begin{rem}
Suppose that $M$ is a Yetter-Drinfeld module over a Hopf algebra $H$ with bijective antipode and denote the braiding on $M$ by $c$. It is easily seen that then $H^{red}(c)$ is a subquotient (i.e. a Hopf algebra quotient of a sub Hopf algebra) of $H$.
\end{rem}

\begin{exmp}
\emph{$H^{red}(c)$ for braidings of group type.} Let $G$ be a group and $M$ a finite dimensional Yetter-Drinfeld module of $G$. Define
\[ C:= \{ g\in G| \exists m\in M:\delta(m)=g\otimes m\}\]
and let $H$ be the subgroup of $G$ generated by $C$. Moreover set
\[   N:= \{ g\in H| \forall m\in M: gm=m\}.\]
Obviously $M$ becomes a Yetter-Drinfeld module over the subquotient $H/N$. It is easy to show that $k(H/N)$ is the reduced FRT construction.
\end{exmp}

% **** Kriterien fuer punktiertheit ****
\section{When is the reduced FRT-construction pointed?}
\label{sect_redpointed}
% **** Bialgebra version ****
The answer to this question was given by Radford in the case of the FRT-bialgebra. In our notation his result reads:
\begin{thm}[\cite{Radford_pointed}, Theorem 3]
\label{thm_bialgpointed}
Let $(M,c)$ be a finite dimensional braided vectorspace. The following are equivalent:
\begin{enumerate}
\item $A^{red}(c)$ is pointed
\item There is a flag of subvectorspaces $0 = M_0\subset M_1\subset\ldots\subset M_r=M$ such that for all $1\leq i\leq r$ $\dim M_i =i$ and
\[c(M_i\otimes M)\subset M\otimes M_i\]
\item There is a flag of $A^{red}(c)$ left subcomodules $0 = M_0\subset M_1\subset\ldots\subset M_r=M$ such that for all $1\leq i\leq r$ $\dim M_i =i$.
\end{enumerate}
\end{thm}

% **** Hopf algebra version ****
We will show now that if $(M,c)$ is rigid we have a similar statement for the reduced FRT Hopf algebra $H^{red}(c)$.
\begin{lem}
\label{lem_cflat}
Let $(M,c)$ be a braided vectorspace and $N\subset M$ a subspace such that
\[ c(N\otimes M)\subset M\otimes N \]
Then we have
\[ c^\flat(M^*\otimes N)\subset N\otimes M^*. \]
\end{lem}
\begin{pf}
This follows directly from the definition of $c^\flat$.
\qed\end{pf}
We will need the following well known statement 
\begin{prop}
\label{prop_leftdual}
Let $H$ be a Hopf algebra and $M$ a left $H$-comodule. Then $M^*$ with the coaction defined by the equation
\[ \forall \varphi\in M^*,m\in M: \varphi\sm1\varphi\s0(m) = S^{-1}(m\sm1) \varphi(m\s0) \]
together with the maps $\ev,\db$ forms a left dual of $M$ in the categorical sense (see e.g. \cite{Kassel}, XIV.2). In particular $c^\flat$ is the inverse of the braiding $c_{M,M^*}$ between $M$ and $M^*$.
\end{prop}
\begin{pf}
It is easy to check that $\ev$ and $\db$ are $H$-colinear. The proof that $c^\flat = c_{M,M^*}^{-1}$ is done as in \cite[XIV.3.1]{Kassel}.
\qed\end{pf}
The following lemma is already used in \cite{Radford_pointed}. We include a proof for completeness.
\begin{lem}
\label{lem_pointedsubcoalg}
Let $H$ be a bialgebra generated (as an algebra) by a subcoalgebra $C\subset H$. If $C$ is pointed then so is $H$.\newline
In this case the coradical of $H$ is generated by the coradical of $C$ as an algebra.
\end{lem}
\begin{pf}
Let $(C_n)_{n\geq 0}$ be the coradical filtration of $C$. Denote by $D_0$ the subalgebra of $H$ generated by $C_0$. As $C$ is pointed we have $D_0\subset kG(H)$. Consider the subsets
\[ D_n := \wedge^n D_0 \:\:\forall n\geq 0.\]
As $D_0$ is a subbialgebra of $H$ the $D_n$ define a bialgebra filtration of the bialgebra
\[ D:= \cup_{n\geq 0} D_n. \]
Now because $C_0\subset D_0$ we have $C_n\subset D_n$ for all $n\geq 0$ and then $C\subset D$. This means $D=H$ and the $D_n$ define a bialgebra filtration of $H$. We find
\[ kG(H)\subset \operatorname{Corad} H \subset D_0\subset kG(H)\]
saying that $H$ is pointed.
\qed\end{pf}

% **** The 1st theorem in Hopf algebra version ***
\begin{thm}
\label{thm_hopfalgpointed}
Let $(M,c)$ be a rigid braided vectorspace. The following are equivalent:
\begin{enumerate}
\item $H^{red}(c)$ is pointed.
\item There is a pointed Hopf algebra $H$ having $M$ as a Yetter-Drinfeld module such that the induced braiding is $c$.
\item There is a flag of left $H^{red}(c)$ subcomodules $0=M_0\subset M_1\subset\ldots\subset M_r = M$ such that for all $1\leq i\leq n$ $\dim M_i = i$.
\item There is a flag of subvectorspaces $0=M_0\subset M_1\subset\ldots\subset M_r = M$ such that for all $1\leq i\leq n$ $\dim M_i = i$ and
\[ c(M_i\otimes M)\subset M\otimes M_i.\]
\end{enumerate}
\end{thm}

\begin{pf}
It is clear that the first item implies the second. If $H$ is as in $(2)$ (e.g. $H=H^{red}(c)$) we find a series of subcomodules 
\[ 0=M_0\subset M_1\subset\ldots\subset M_r=M \]
with $\dim M_i = i$ for all $1\leq i\leq r$. Together with the definition of the braiding for Yetter-Drinfeld modules this gives us that $(1)$ implies $(3)$ (and hence also $(4)$) and that $(2)$ implies $(4)$. We still have to show that $(4)$ implies $(1)$.\newline
So now assume that $(4)$ holds.
For all $1\leq i\leq r$ choose $m_{r+1-i}\in M_i\setminus M_{i-1}$ arbitrarily (thus $m_i,\ldots,m_r\in M_i$). This defines a basis of $M$ and let $m^1,\ldots,m^r$ be the dual basis. We find elements $t_{ij}\in H(c), 1\leq i,j\leq r$ satisfying
\[ \Delta(t_{ij}) = \sum\limits_{l=1}^r t_{il}\otimes t_{lj}\:\:\mbox{and}\:\: \delta(m_i) = \sum\limits_{l=1}^r t_{il}\otimes m_l. \]
Using the definition of the coaction of $M^*$ we see that
\[ \delta(m^i) = \sum\limits_{l=1}^r S^{-1}(t_{li})\otimes m^l. \]
Now define
\[ J := \mbox{k-span}\{t_{ij} | 1\leq j < i\leq r \}. \]
$J$ is a coideal of $H(c)$ and we will show $J\subset J_{H(c)}$. 
For all $1\leq i,k \leq r$ we have
\[ c(m_k\otimes m_i) \in c(M_k\otimes M)\subset M\otimes M_k\]
and on the other hand
\[ c(m_k\otimes m_i) = \sum\limits_{l=1}^r t_{kl}m_i\otimes m_l 
= \sum\limits_{j,l=1}^r r(t_{ij},t_{kl}) m_j\otimes m_l\]
This implies $r(t_{ij},t_{kl})=0$ for $l<k$. Moreover because of $c(M_i\otimes M)\subset M\otimes M_i$ we have by Lemma \ref{lem_cflat} that $c^\flat(M^*\otimes M_i)\subset M_i\otimes M^*$ and thus by Proposition \ref{prop_leftdual} that
$$ c_{M,M^*}(M_i\otimes M^*)\subset M^*\otimes M_i$$
In the same manner as above we obtain $r(S^{-1}(t_{ij}),t_{kl})=0$ for $l< k$.Now by Theorem \ref{thm_FRTHA} we have that the algebra $H(c)$ is generated by the $t_{ij}$ and the $S(t_{ij})$. Thus it is also generated by the $S^{-1}(t_{ij})$ and the $t_{ij}$. Lemma \ref{lem_radical_gen} allows us to conclude that $r(H\otimes J)=0$ and thus $J\subset J_{H(c)}$. As $J_{H(c)}$ is stable under $S$ we obtain $J+S(J)\subset J_{H(c)}$.\newline
To see that $H^{red}(c)$ is pointed it suffices to show that the coalgebra $C$ spanned by the images of $\overline{t_{ij}},S(\overline{t_{kl}})$ of $t_{ij},S(t_{kl})$ in $H^{red}(c)$ is pointed (Lemma \ref{lem_pointedsubcoalg}). For this define subsets $C_n,n\geq 0$ by
\[ C_n := \operatorname{k-span}\{\overline{t_{ij}},S(\overline{t_{ij}})|1\leq i\leq j\leq i+n\leq r \}\subset H^{red}(c)\]
Because $J+S(J)\subset J_{H(c)}$ we find that the $C_n$ define a coalgebra filtration of $C$. Thus 
\[ \operatorname{Corad} C\subset C_0 = \operatorname{k-span}\{\overline{t_{ii}},S(\overline{t_{ii}})|1\leq i\leq r\}.\]
As $\Delta(\overline{t_{ii}})= \overline{t_{ii}}\otimes \overline{t_{ii}}$ and $\Delta(S(\overline{t_{ii}})) = S(\overline{t_{ii}})\otimes S(\overline{t_{ii}})$ in $H^{red}(c)$ we find that $C$ is pointed.
\qed\end{pf}
For future usage we remark that the coradical of $H^{red}(c)$ is generated by the elements $\overline{t_{ii}},S(\overline{t_{ii}})$.

% **** Triangulaere Verzopfungen ****
\section{The reduced FRT construction for triangular braidings}
\label{sect_triang}
In this section we want to refine our knowledge of the reduced FRT bialgebra and FRT Hopf algebra for a special class of braidings considered in \cite{Ich_PBW}.

% **** Definition triangulaer ****
\begin{defn}
\label{defn_triangular}
Let $(M,c)$ be a finite dimensional vectorspace with a totally ordered basis $X$. $(M,c)$ will be called \emph{left triangular} (with respect to the basis $X$) if for all $x,y,z\in X$ with $z\gx y$ there exist $\gamma_{x,y}\in k\setminus \{0\}$ and $v_{x,y,z}\in M$ such that
\[ c(x\otimes y) = \gamma_{x,y} y\otimes x + \sum\limits_{z\gx y} z\otimes v_{x,y,z} \mtxt{for all} x,y\in X \]
$(M,c)$ will be called \emph{right triangular} (with respect to the basis $X$) if for all $x,y,z\in X$ with $z\gx x$ there exist $\beta_{x,y}\in k\setminus \{0\}$ and $w_{x,y,z}\in M$ such that
\[ c(x\otimes y) = \beta_{x,y} y\otimes x + \sum\limits_{z\gx x} w_{x,y,z}\otimes z \mtxt{for all} x,y\in X \]
\end{defn}

% **** Proposition c^{-1}<-->c ****
\begin{prop}
\label{prop_c_c-1}
Let $(M,c)$ be a finite dimensional braided vectorspace.
\begin{enumerate}
\item $c$ is left triangular if and only if $c^{-1}$ is right triangular.
\item $c$ is left triangular if and only if $\tau c\tau$ is right triangular.
\end{enumerate}
\end{prop}
\begin{pf}
$(2)$ is trivial. Thus for $(1)$ it suffices to show the if-part. Assume $c^{-1}$ is right triangular and adopt the notation from the definition. Define $M_{\gx x} := \operatorname{k-span}\{z\in X|z\gx x\}$. We see from the definition that
\[ c(y\otimes x) = \beta_{xy}^{-1} x\otimes y + c\left(\sum\limits_{z\gx x} w_{x,y,z}\otimes z\right) \in \beta_{xy}^{-1} x\otimes y + c(M\otimes M_{\gx x})\]
It is now easy to show by downward induction on $x$ (along the order on $X$)
\[ c(y\otimes x)\in \beta_{xy}^{-1}x\otimes y + M_{\gx x}\otimes M. \]
Thus $c$ is left triangular.
\qed\end{pf}

% **** Die sind starr ****
The first important property of triangular braidings in our context is that they are rigid. Thus the notion of the (reduced) FRT construction makes sense.
\begin{lem}
\label{lem_rigid}
Let $(M,c)$ be a (left or right) triangular braided vectorspace. Then $(M,c)$ and $(M,c^{-1})$ are rigid.
\end{lem}
\begin{pf}
In both cases it suffices to show that $(M,c)$ is rigid. Assume $(M,c)$ is left triangular with respect to the basis $X$. Let $(\varphi_x)_{x\in X}$ denote the dual basis ($\varphi_x(y)=\delta_{xy}$). Then
\begin{eqnarray*}
c^\flat\tau(x\otimes\varphi_y) &=& c^\flat(\varphi_y\otimes x)\\
&=& \sum\limits_{z\in X} (\varphi_y\otimes M)c(x\otimes z)\otimes \varphi_z \\
&=& \gamma_{x,y} x\otimes\varphi_y + \sum\limits_{z\in X,z'\gx z} \varphi_y(z') v_{x,z,z'}\otimes \varphi_z\\
&\in& \gamma_{x,y} x\otimes\varphi_y + \sum\limits_{z\lx y} v_{x,z,z'}\otimes \varphi_z
\end{eqnarray*}
This means that the map $c^\flat\tau$ has upper triangular representing matrix with respect to the basis
\[ x_1\otimes\varphi_{x_1},\ldots,x_r\otimes\varphi_{x_1},x_1\otimes\varphi_{x_2},\ldots,x_r\otimes\varphi_{x_2},\ldots,x_1\otimes\varphi_{x_r},\ldots,x_r\otimes\varphi_{x_r}\]
(where we assumed that the elements of $X$ are $x_1\lx x_2\lx\ldots\lx x_r$). The diagonal entries are $\gamma_{x,y}\neq 0$ and thus the matrix is invertible. This shows that $c^\flat$ is an isomorphism. A similar proof works for right triangular braidings.
\qed\end{pf}

These results together with the Theorems \ref{thm_bialgpointed}, \ref{thm_hopfalgpointed} already show that the reduced FRT-construction $A^{red}(c)$ (resp. $H^{red}(c)$) is pointed if $c$ is right triangular. We will describe these reduced FRT-constructions more exactly. To formulate this we introduce

% **** Definition Eigenspaces ****
\begin{defn}
Let $G$ be an abelian monoid and $M$ a $G$-module. We say \emph{$G$ acts diagonally on $M$} if $M$ is the direct sum of simultaneous eigenspaces under the action of $G$, this means:
\[ M = \bigoplus\limits_{\chi\in\hat{G}}\{m\in M| \forall g\in G: gm = \chi(g)m\}.\]
\end{defn}
% **** Die G-Kompatible Fahne ****
We will need the following technical facts:
\begin{prop}
\label{prop_flag}
Let $H$ be a pointed bialgebra with abelian coradical such that for all $g\in G:=G(H)$ the map $H\rightarrow H,h\mapsto hg$ is injective. Let $M\in{}^H_H\mYD$ be finite-dimensional and assume that $G$ acts diagonally on $M$.
Then there is a series of $H$-subcomodules and $G$-submodules
\[ 0 = M_0\subset M_1\subset\ldots\subset M_r = M\]
such that $\dim M_i = i$ for all $1\leq i\leq r$.
\end{prop}
\begin{pf}
Consider modules $N\in {}^H\mM\cap {}_G\mM$ that have an eigenspace decomposition as in the lemma and satisfy the following compatibilty condition:
\[ (gn)_\sm1 g\otimes (gn)_\s0 = gn_\sm1\otimes gn_\s0 \:\forall g\in G,n\in N.\]
It suffices to show that every such module $N$ contains a one dimensional $H$-subcomodule that is also a $G$-submodule (note that the objects considered in the lemma are of this type).\newline
So pick a simple subcomodule of $N$. As $H$ is pointed this is spanned by some $n_0\in N$ and we find $g\in G$ such that $\delta(n_0)=g\otimes n_0$. Consider the subcomodule
\[ 0\neq X := \{n\in N|\delta(n)=g\otimes n\}\subset N. \]
This is a $G$-submodule as for $n\in X$ and $h\in G$ we have by the compatibility condition (and because $G$ is abelian)
\[ (hn)_\sm1 h\otimes (hn)_\s0 = hn_\sm1\otimes hn_\s0 = hg\otimes hn = gh\otimes hn. \]
Now right multiplication with $h$ is injective by assumption and we obtain $\delta(hn) = g\otimes hn$ showing that $X$ is indeed a $G$-submodule. A lemma from linear algebra tells us that because $N$ is the direct sum of eigenspaces under the action of $G$, so is $X$. We find an element $n\in N$ that is an eigenvector under the action of $G$. Then $kn$ is a one dimensional $H$-subcomodule and $G$-submodule.
\qed\end{pf}
% **** Prop Komodulstruktur ****
\begin{prop}
\label{prop_YDbasis}
Let $H$, $M$ be as in \ref{prop_flag}. Then we can find a basis $m_1,\ldots,m_r$ of $M$ made up of eigenvectors under the action of $G(H)$ and elements $c_{ij}\in H,1\leq i\leq j\leq r$ such that
\[ \delta(m_i) = \sum\limits_{l=i}^r c_{il}\otimes m_l,\:\:\Delta(c_{ij})=\sum\limits_{l=i}^j c_{il}\otimes c_{lj},\:\:\eps(c_{ij}) = \delta_{ij}\]
\end{prop}
\begin{pf}
Take the series of comodules from \ref{prop_flag}. For all $1\leq i\leq r$ we can choose an eigenvector $m_{r+1-i}\in M_i\setminus M_{i-1}$ (thus $M_i=\operatorname{k-span}\{m_i,\ldots,m_r\}$). Now we can find elements $c_{ij}\in H,1\leq i\leq j\leq r$ such that
\[ \delta(m_i) = \sum\limits_{l=i}^r c_{il}\otimes m_l. \]
The formulas for the comultiplication and the counit follow from the axioms of the comodules.
\qed\end{pf}
% **** Hopf algebra version ****
\begin{thm}
\label{thm_hopfalg}
Let $(M,c)$ be a rigid braided vectorspace. The following are equivalent:
\begin{enumerate}
\item $c$ is right triangular
\item $H^{red}(c)$ is pointed with abelian coradical and $G(H^{red}(c))$ acts diagonally on $M$.
\item There is a pointed Hopf algebra $H$ with abelian coradical having $M$ as a Yetter-Drinfeld module such that the induced braiding is $c$ and $G(H)$ acts diagonally on $M$.
\end{enumerate}
\end{thm}
\begin{pf}
Of course $(2)$ implies $(3)$. Assume $c$ is right triangular with respect to the basis $m_1,\ldots,m_r$ ordered by $m_1\lx \ldots\lx m_r$. let $M_i := \operatorname{k-span}\{m_i,\ldots,m_r\}$. Then we have of course
\[ c(M_i\otimes M)\subset M\otimes M_i\:\:\forall 1\leq i\leq r. \]
Theorem \ref{thm_hopfalgpointed} tells us that $H^{red}(c)$ is pointed. We adopt the notation from the proof of $(4)\Rightarrow (1)$ there. Then we obtain using the right triangularity of $c$:
\[ \sum\limits_{l=i}^r \overline{t_{il}}m_j\otimes m_l = c(m_i\otimes m_j)\in \alpha_{ij} m_j\otimes m_i + M\otimes\sum\limits_{l=i+1}^r k m_l. \]
This means $\overline{t_{ii}}m_j = \alpha_{ij}m_j$ for all $1\leq i,j\leq r$. As the $\overline{t_{ii}}$ and their inverses generate the coradical of $H^{red}(c)$ as an algebra we get that $G(H^{red}(c))$ acts diagonally on $M$.\newline
We are left to show that the $G(H^{red}(c))$ is abelian. Let $g,h\in G(H^{red}(c))$, thus $g$ and $h$ act diagonally on $M$. Then $gh-hg$ acts as $0$ on $M$, saying $k(gh-hg)$ is a coideal annihilating $M$. In the same way $g^{-1}h^{-1}-h^{-1}g^{-1}$ annihilates $M$ and thus $k(gh-hg)$ annihilates $M^*$. As $H^{red}(c)$ is $M,M^*$-reduced, we get that $gh=hg$.\\[0,2cm]
Now assume we are given a Hopf algebra as in $(3)$. Let $m_1,\ldots,m_r$ be the basis of $M$ from Proposition \ref{prop_YDbasis} and also take the $c_{ij}$ from there. Then we have
\[ c(m_i\otimes m_j) \in c_{ii}m_j\otimes m_i + M\otimes\sum\limits_{l>i}km_l\]
as $c_{ii}\in G(H)$ we have $c_{ii}m_j\in km_j\setminus\{0\}$ and thus $c$ is right triangular.
\qed\end{pf}

\begin{cor}
Assume that $k$ is algebraically closed. Let $G$ be an abelian group and $M$ a finite-dimensional Yetter-Drinfeld module over $G$. Then the induced braiding is left triangular. If it is also right triangular, then it is of diagonal type, i.e. there is a basis $X$ of $M$ and there are nonzero scalars $q_{xy}\in k$ for all $x,y\in X$ such that the braiding is given by
\[ c(x\otimes y)=q_{xy}y\otimes x.\]
In particular there are braidings which are left but not right triangular. For example take $k=\C$, $G=\langle g \rangle\isom \Z$ and let $g$ act on a two-dimensional vectorspace $M$ by a Jordan block ($\lambda\in \C\setminus\{0\}$)
\[ \Matrix{cc}{\lambda & 1\\
0&\lambda\\}.\]
The coaction is given by $\delta(m) := g\otimes m$ for all $m\in M$.
\end{cor}
\begin{pf}
For all $g\in G$ let $M_g := \{m\in M|\delta(m) = g\otimes m\}$. Then the $M_g$ are $G$-submodules of $M$. Since every simple submodule of a finite-dimensional $G$-module is one-dimensional we see that each $M_g$ has a flag of invariant subspaces. So for each $g\in G$ we find a basis $m^g_1,\ldots,m^g_{r_g}$ of $M_g$ such that for all $h\in G$
\[h\cdot m^g_i\in km^g_i\oplus\ldots\oplus km^g_{r_g}.\]
Now by concatenating these bases and ordering each according to the indices we obtain a totally ordered basis such that $c$ is triangular.\\
Now assume that the braiding is right triangular. By passing to a quotient of $G$ we may assume that $gm=m$ for all $m\in M$ implies $g=1$. We show that $kG$ is $M,M^*$ reduced. Let $J\subset kG$ be the sum of all coideals annihilating $M$ and $M^*$. It is easy to see that $J$ is a Hopf ideal. Then $kG/J\isom kH$ is  a group algebra of some group $H$ and $M$ is again a Yetter-Drinfeld module over $kH$. So we get an induced epimorphism of groups $\pi:G\rightarrow H$. For $g\in\ker\pi$ we have $gm=m$ for all $m\in M$ and thus $g=1$. This means that $\pi$ is injective, hence $J=0$.\\
By the universal property $H^{red}(c)$ is a sub Hopf algebra of $kG$ and hence a group algebra as well. But $G(H^{red}(c))$ acts diagonally on $M$, which means that the braiding is indeed of diagonal type.
\qed\end{pf}
% **** Bialgebra version ****
% There is also a version for the reduced FRT-bialgebra. The proof is very similar to the one for theorem \ref{thm_hopfalg} and will be omitted.
% \begin{thm}
% \label{thm_bialg}
% Let $(M,c)$ be a finite dimensional braided vectorspace. The following are equivalent:
% \begin{enumerate}
% \item $c$ is right triangular
% \item The following three conditions hold
% \begin{itemize}
% \item $A^{red}(c)$ is pointed with abelian coradical 
% \item $G(A)$ acts diagonally on $M$
% \item for all $g\in G$ the map $\tilde{g}:A\rightarrow A, h\mapsto hg$ is injective.
% \end{itemize}
% \end{enumerate}
% \end{thm}

\section{$H^{red}(c)$ for finite dimensional $U_q(\mathfrak{g})$-modules}
\label{sect_uqg}
In this section we will determine $H^{red}(c)$ for braidings induced by finite dimensional $U_q(\mathfrak{g})$-modules. Assume that $k$ is an algebraically closed field of characteristic zero. In this section we use details from \cite[sections 2,3 and 4]{Ich_uqg} and several statements on the representation theory of $U_q(\mathfrak{g})$; a good reference for the latter is \cite{Jantzen}.\par
We will need the following proposition for Radford biproducts (the notations for Radford biproducts are taken from \cite[section 2]{Ich_uqg} or, equivalently, from \cite{AS5}).
\begin{prop}
\label{prop_biprod}
Assume that $\psi:A\rightarrow A'$ is a morphism between Hopf algebras $A,A'$ with bijective antipodes. Let $H\subset A,H'\subset A'$ be sub Hopf algebras with Hopf algebra projections $p,p'$ such that the following diagram
\begin{diagram}
A&\rTo^\psi&A'\\
\dTo<p&&\dTo>p'\\
H&\rTo^{\psi|H}&H'
\end{diagram}
commutes (and is well defined, i.e. $\psi(H)\subset H'$). Let $R:=A^{\co p},R':=A'^{\co p'}$. Then $\psi(R)\subset R'$ and the diagram
\begin{diagram}
A&\rTo^\psi&A'\\
\dTo<\isom&&\dTo\isom\\
R\#H&\rTo_{\psi|R\#\psi|H}&R'\#H'
\end{diagram}
commutes, where the vertical isomorphisms are given by
\[ A\rightarrow R\#H, a\mapsto a\s1 S_Hp(a\s2)\#p(a\s3)\]
and the corresponding map for $A'$.
\end{prop}
\begin{pf}
The vertical isomorphisms are those from Radford's theorem on Hopf algebras with a projection \cite{Radford_HAproj}. The rest of the proposition is just a computation. \qed
\end{pf}

Now let $M$ be a finite dimensional $U_q(\mathfrak{g})$-module with braiding $c=c^f$. We will use the notations from \cite[sections 3, 4]{Ich_uqg}. Define
\[ P:= \{\alpha\in\Pi| E_\alpha M \neq 0\}, \]
\[ W:= \{\lambda\in\Lambda | M_\lambda \neq 0\}.\]
Let $\hat{U}=\Nichols(V)\#kG$ be the Hopf algebra defined in \cite[section 4]{Ich_uqg}, where $V$ is a completely reducible Yetter-Drinfeld module over the abelian group $G$. Let $H$ be the subalgebra of $\hat{U}$ generated by the $\hat{F}_\alpha,\alpha\in P$ and the $K_\lambda^{\pm1},\lambda\in W$. Denote by $\tilde{V}$ the subspace of $V$ generated by the $\hat{F}_\alpha,\alpha\in P$; this is again a Yetter-Drinfeld module over $G$. Furthermore define $\tilde{G}:=G(H)$ and
\[ N := \{g\in\tilde{G}|\forall m\in M: gm=m\},\]
\[ J := \operatorname{k-span}\{gn-g|g\in G,n\in N\setminus\{1\}\}.\]
The aim of this section is to prove the following theorem:
\begin{thm}
The reduced FRT construction of $(M,c^f)$ is given by
\[ H^{red}(c^f) \isom \Nichols(\tilde{V})\#k(\tilde{G}/N).\]
\end{thm}

\begin{pf}
First observe that the Yetter-Drinfeld module $V$ over $G$ can be restricted to a Yetter-Drinfeld module over $\tilde{G}$ because $K_\alpha\in \tilde{G}$ for all $\alpha\in P$: For $\alpha\in P$ we find $\lambda\in W,m\in M_\lambda$ with $0\neq\hat{F}_\alpha m\in M_{\lambda-\alpha}$. By the definition of $\hat{U}$ and of the coaction on $M$ it follows that $L_\lambda,L_{\lambda-\alpha}\in W\subset \tilde{G}$. Hence also $K_\alpha = L_{\lambda-\alpha}L_\lambda^{-1}\in\tilde{G}$.\par
Next we show that $N$ acts trivially on $V$. Let $g\in N,\alpha\in P$; then there is an $m\in M$ with $\hat{F}_\alpha m\neq 0$ and there is an $\eps\in k$ such that $g\cdot\hat{F}_\alpha = \eps \hat{F}_\alpha$. Then
\[\eps \hat{F}_\alpha m = (g\cdot\hat{F_\alpha})m = g \hat{F}_\alpha g^{-1} m = \hat{F}_\alpha m\]
implies $\eps = 1$ and thus $g\cdot\hat{F}_\alpha=\hat{F}_\alpha$.
This means that $V$ is can be turned into a Yetter-Drinfeld module over $\tilde{G}/N$ using the canonical projection $\tilde{G}\rightarrow \tilde{G}/N$. Note that $k(\tilde{G}/N)$ is the reduced FRT construction for $\tilde{V}$. Now we can form $\tilde{H} := \Nichols(\tilde{V})\#k(\tilde{G}/N)$. We have a canonical projection $H=\Nichols(\tilde{V})\#k\tilde{G}\rightarrow \tilde{H}$.\par
Now observe that the $\hat{U}$ coaction on $M$ can be restricted to a $H$ coaction. As $N$ acts trivially on $M$ by definition, we can turn $M$ into a Yetter-Drinfeld module over $\tilde{H}$ using the canonical projection. We obtain then a commutative diagram of Hopf algebra projections
\begin{diagram}
H(c)&\rOnto^\varphi&\tilde{H}&=\Nichols(\tilde{V})\#k(\tilde{G}/N)\\
\dOnto<\pi&\ldOnto_\psi&&\\
H^{red}(c)&&&\\
\end{diagram}
where $\varphi$ is given by the universal property of $H(c)$ and $\pi$ is the canonical projection. Both maps are compatible with action and coaction. To show that we have a factorization $\psi$ we show $\ker\varphi\subset\ker\pi$: Let $x\in\ker\varphi$. Then $xM=\varphi(x)M=0,xM^*=\varphi(x)M^*=0$. This implies that $x\in J_{H(c)}=\ker\pi$.\par
To show that $\psi$ is injective we will first show that all occurring maps are graded. Recall from \cite[section 6]{Ich_uqg} that $M$ has a $\N$-grading such that the structure maps of the $\hat{U}$ Yetter-Drinfeld module structure are graded ($\hat{U}$ is coradically graded by construction). This grading turns the $\tilde{H}$ action and coaction into graded maps and shows that the braiding $c$ is graded. Thus $H(c)$ and $H^{red}(c)$ have $\Z$ gradings such that the projection $\pi$, the actions and the coactions are graded (This can easily be seen in the construction of $H(c)$ given in \cite{Sbg_cqt}: Start with a homogenous basis $m_1,\ldots,m_r$ of $M$ and grade $H(c)$ by giving the generator $T_{ij}$ the degree $deg(m_i)-deg(m_j)$). Using the compatibility condition between $\varphi$ and the $H(c)$ resp. $\tilde{H}$ coactions it is easy to see that also $\varphi$ is a graded map and hence $H(c)$ and $H^{red}(c)$ are actually $\N$ graded. By construction also the map $\psi$ is graded.\par
Now both $\tilde{H}$ and $H^{red}(c)$ are graded Hopf algebras, hence admit Hopf algebra projections onto the zeroth components. As $\psi$ is a graded map we can apply Proposition \ref{prop_biprod} to our situation and see that in order to show that $\psi$ is injective it suffices to show that $\psi|k(\tilde{G}/N)$ and $\psi|\Nichols(\tilde{V})$ are injective.\par
First show that $\psi|k(\tilde{G}/N)$ is injective: Let $\bar{x},\bar{y}\in \tilde{G}/N$ such that $\psi(\bar{x})=\psi(\bar{y})$. This means $xm=ym$ for all $m\in M$ and thus $xy^{-1}\in N$. Hence $\bar{x}=\bar{y}$, showing that $\psi|\tilde{G}/N$ is injective. The claim follows by linear algebra.\par
On the other hand, let $I$ be the kernel of $\psi|\Nichols(\tilde{V})$. As this is a graded morphism of algebras and coalgebras, $I$ is a coideal and an ideal generated by homogeneous elements. By the characterization of Nichols algebras from \cite{AS5}, $I=0$ if $I\cap \tilde{V}=0$ (i.e. $I$ is generated by elements of degree $\geq 2$).\par
So assume we have $x\in I\cap \tilde{V}$ and write $x=\sum\limits_{\alpha\in P}r_\alpha \hat{F}_\alpha$ for scalars $r_\alpha\in k$. Then, as $I\subset\ker\psi$, $xM=0$. The weight-space grading of the modules $M$ yields that for all $\alpha\in P$
\[ r_\alpha \hat{F}_\alpha M=0.\]
As $\hat{F}_\alpha M\neq 0$ by construction of $P$ we get $r_\alpha=0$ for all $\alpha\in P$ and thus $x=0$.
\qed\end{pf}

\begin{rem}
The set $P$ is a union of connected components of the Coxeter graph of $\mathfrak{g}$.\\
In particular if $\mathfrak{g}$ is simple, we have $\tilde{V}=V$ and thus $H^{red}(c)$ is obtained from $H$ just by dividing out the ideal generated by the set
\[ \{g-h| g,h\in G,\forall m\in M gm=hm\}.\]
In the general case we obtain that $H^{red}(c)$ may be viewed as the ``nonpositive part of a quantized enveloping algebra of $\hat{\mathfrak{g}}$'' (where $\hat{\mathfrak{g}}$ is the sub Lie algebra of $\mathfrak{g}$ generated by the $E_\alpha,H_\alpha,F_\alpha, \alpha\in P$) in the sense that $H^{red}(c)$ is a biproduct of the negative part $U_q^-(\hat{\mathfrak{g}})$ with a finitely generated  abelian group of finite rank. The group algebra of this group is exactly the reduced FRT construction of $\tilde{V})$.
\end{rem}
\begin{pf}
The second part of follows from the proof of the theorem above. We show only that $P$ is a union of connected components of the Coxeter graph, i.e. if $\alpha\in P$ and $\beta\in\Pi$ with $(\alpha,\beta)<0$ then also $\beta\in P$. \\
So assume we have $\alpha\in P,\beta\in\Pi$ such that $(\alpha,\beta)<0$. Thus we have $E_\alpha M\neq 0$ and we will show $E_\beta M\neq 0$. Let $\lambda\in\Lambda,m\in M_\lambda$ with $E_\alpha m\neq 0$. If $(\lambda,\beta)=0$ replace $m$ by $E_\alpha m$ and $\lambda$ by $\lambda + \alpha$. Hence $0\neq m\in M_\lambda$ and $(\lambda,\beta)\neq 0$. Let $U_q(\mathfrak{sl}_2)_\beta$ be the subalgebra of $U_q(\mathfrak{g})$ generated by $E_\beta,F_\beta,K_\beta$ and $K_\beta^{-1}$; it is isomorphic to $U_{q^{(\beta,\beta)}}(\mathfrak{sl}_2)$ as a Hopf algebra. Consider the $U_q(\mathfrak{sl}_2)_\beta$ submodule $N$ of $M$ generated by $m$. If $E_\beta m\neq 0$ we have $\beta \in P$ and the proof is done. So assume $E_\beta m=0$, hence $m$ is a highest weight vector for the  $U_q(\mathfrak{sl}_2)_\beta$-module $N$. As $(\lambda,\beta)\neq 0$, $N$ is not one-dimensional. Thus we have $E_\beta(F_\beta m)\neq 0$ implying $\beta \in P$.\qed
\end{pf}

\begin{rem}
It is an open question if there is a combinatorial description of those triangular braidings for which the reduced FRT construction is generated by grouplike and skew-primitive elements. 
\end{rem}

\bibliography{../../tex/promotion.bib}

\end{document}